\documentclass{article}
\usepackage{amsmath,amsthm,amsopn,amstext,amscd,amsfonts,amssymb}
\usepackage{natbib}
\usepackage{verbatim}

\def \tr {\mathop{\rm tr}\nolimits}
\def \etr {\mathop{\rm etr}\nolimits}
\def \diag {\mathop{\rm diag}\nolimits}

\def \R{\mathop{\rm Re}\nolimits}

\renewenvironment{abstract}
                 {\vspace{6pt}
                  \begin{center}
                  \begin{minipage}{5in}
                  \centerline{\textbf{Abstract}}
                  \noindent\ignorespaces
                 }
                 {\end{minipage}\end{center}}
\newtheorem{thm}{\textbf{Theorem}}[section]
\newtheorem{cor}{\textbf{Corollary}}[section]
\newtheorem{lem}{\textbf{Lemma}}[section]
\theoremstyle{definition}

\newtheorem{rem}{\textbf{Remark}}[section]

\title{\huge \textbf{A note about zonal polynomials}}
\author{
  \textbf{Jos\'e A. D\'{\i}az-Garc\'{\i}a} \thanks{Corresponding author\newline
    }\\
  Department of Statistics and Computation \\
  25350 Buenavista, Saltillo, Coahuila, Mexico \\
  E-mail: jadiaz@uaaan.mx \\[2ex]
  \textbf{Ram\'on Guti\'errez J\'aimez} \\
  Department of Statistics and O.R. \\
  University of Granada \\
  Granada 18071, Spain \\
  E-mail: rgjaimez@ugr.es\\
}
\date{}
\begin{document}
\maketitle
\begin{abstract}
In this paper, we study some properties of multivariate gamma function and zonal
polynomials.
\end{abstract}

\section{\large{Introduction}}

Zonal polynomials are of undeniable importance, in both theory and practice.
Following the algorithms proposed by \citet{KE06}, zonal polynomials are increasingly
being used in various areas of knowledge. Undoubtedly, the initial studies by
\citet{j:60}, \citet{j:61a}, \citet{j:61b}, \citet{c:63} and \citet{j:64}, among
others, laid the foundations in this field. Subsequently, books by \citet{f:85} and
\citet{t:84}, among others, compiled many of these early results and proposed new
theoretical considerations and many practical applications. Nevertheless, the book by
\citet{mh:82}, above all others, marks a watershed in these studies and has had an
undeniable impact on recent generations of mathematicians and statisticians working
in the field of multivariate analysis, see \citet{w:84}. Virtually all recent studies
that bear upon the question of zonal polynomials have cited Muirhead's book (1982).
In particular, \citet{R:97}, calculated the expectation of zonal and invariate
polynomials, making use of various results published in \citet{mh:82}. Unfortunately,
the conclusions drawn by \citet{R:97} are incorrect, and this is because both Lemma
7.2.12 and the proof of Theorem 7.2.13 in \citet{mh:82} are incorrect. Due to the
undeniable importance and impact of Muirhead's book, and its influence on current and
future studies, the present text proposes corrections to the above-mentioned lemma
and theorem (see Section 3). Prior to this, in Section 2 we propose an expression and
alternative proof of the multivariate gamma function.

\section{\large{Preliminary results}}\label{sec1}

The Pochhammer symbol is defined as
\begin{equation}\label{eq1}
    (x)_{q}= x(x+1)\cdots (x + q -1) = \prod_{i = 1}^{q}(x + i - 1) = \prod_{i = 1}^{q}(x + q -
  i) = \frac{\Gamma[x + q]}{\Gamma[x]},
\end{equation}
where $\Gamma[\cdot]$ is the gamma function. Also, observe that
\begin{equation}\label{negpch}
    (-x)_{q} = (-1)^{q}(x-q+1)_{q} = \frac{(-1)^{q}\Gamma[x+1]}{\Gamma[x-q+1]}.
\end{equation}
Similarly
$$
  (x)^{(q)} = x(x-1)\cdots (x - q + 1) = \prod_{i = 1}^{q}(x-i+1) = \prod_{i = 1}^{q}(x-q
  +i).
$$

Then for any function $g:\Re \rightarrow \Re$
\begin{equation}\label{prod1}
    \prod_{i = 1}^{q}g(x + i - 1) = \prod_{i = 1}^{q}g(x + q - i),
\end{equation}
and
\begin{equation}\label{prod2}
    \prod_{i = 1}^{q}g(x-i+1) = \prod_{i = 1}^{q}g(x-q +i).
\end{equation}

\begin{lem} \label{lemma1}Let $X$ be a real $m \times m$ positive definite matrix and $\R(a) > (m-1)/2$.
The multivariate gamma function, denoted by $\Gamma_{m}[a]$, is defined to be
$$
  \Gamma_{m}[a] = \int_{A>0} \etr(-A) (\det A)^{a-(m+1)/2}(dA),
$$
where $\etr (\cdot) \equiv \exp\tr(\cdot)$ and the integral is over the space of
positive definite (and hence symmetric) $m \times m$ matrices. Here, $A > 0$ means
that $A$ is a positive definite matrix. Then
\begin{equation}\label{gama1}
    \Gamma_{m}[a] = \pi^{m(m-1)/4}\prod_{i = 1}^{m} \Gamma[a-(i-1)/2],
\end{equation}
\begin{equation}\label{gama2}
    \phantom{\Gamma_{m}[a]} = \pi^{m(m-1)/4}\prod_{i = 1}^{m} \Gamma[a-(m-i)/2].
\end{equation}
\end{lem}
\textit{Proof.} (\ref{gama1}) is given in \citet[Theorem 2.1.12, pp. 62-63]{mh:82},
see also, \citet[Example 1.24, pp. 56-57]{m:97}.

For proof (\ref{gama2}), let $T$ be a real upper-triangular matrix with $t_{ii} > 0$,
$i = 1,\dots,m$ and consider the decomposition $X = TT'$, then
$$
  \tr A = \tr TT' = \tr T'T = \sum_{i\leq j}^{m} t_{ij}^{2},
$$
$$
  \det A = \det TT' = (\det T)^{2} = \prod_{i = 1}^{m}t_{ii}^{2},
$$
and from \citet[Theorem 1.28, p. 56]{m:97}
$$
  (dA) = 2^{m}\prod_{i = 1}^{m}t_{ii}^{i} \bigwedge_{i \leq j}^{m}dt_{ij} = \prod_{i = 1}^{m}\left
  (t_{ii}^{2}\right)^{(i-1)/2} \bigwedge_{i = 1}^{m}dt_{ii}^{2}\bigwedge_{i < j}^{m}dt_{ij}.
$$
Hence,
$$
  \Gamma_{m}[a] = \prod_{i < j}^{m} \left(\int_{-\infty}^{\infty} \exp\left(-t_{ij}^{2}\right)dt_{ij}\right )
  \prod_{i = 1}^{m} \left (\int_{0}^{\infty} \exp\left(-t_{ii}^{2}\right) \left (t_{ii}^{2} \right )^{a-(m-i)/2- 1} dt_{ii}^{2}
  \right).
$$
But
$$
  \int_{-\infty}^{\infty} \exp\left(-t_{ij}^{2}\right)dt_{ij} = \sqrt{\pi}
$$
and
$$
  \int_{0}^{\infty} \exp\left(-t_{ii}^{2}\right) \left (t_{ii}^{2} \right )^{a-(m-i)/2- 1}
  dt_{ii}^{2} = \Gamma[a-(m-i)/2].
$$
From where (\ref{gama2}) follows. Alternatively, (5) is an immediate consequence of
(\ref{prod2}). \qed

In a similar way to expression (\ref{prod2}) it is readily apparent that for $k_{i}$
non negative integers, $i = 1, \dots,q$,
\begin{equation}\label{prod3}
    \prod_{i = 1}^{q}g(x \pm k_{q+1-i} -i+1) = \prod_{i = 1}^{q}g(x \pm k_{i}-q +i).
\end{equation}

\section{\large{Zonal polynomials}}\label{sec2}

In this section we propose the correct version of Lemma 7.2.12, p. 256 and the
correct proof of Theorem 7.2.13, pp. 256-258 in \citet{mh:82}.

\begin{lem} \label{lemma2}If $Z = \diag(z_{1}, \dots, z_{m})$ and $Y =
(y_{ij})$ is an $m \times m$ positive definite matrix then
\begin{eqnarray*}
  C_{\kappa}\left(Y^{-1}Z\right) &=& d_{\kappa}z_{1}^{k_{1}} \cdots z_{m}^{k_{m}} y_{11}^{-(k_{m} -
  k_{m-1})} \det
  \left[
    \begin{array}{cc}
      y_{11} & y_{12} \\
      y_{21} & y_{22}
    \end{array}
  \right ]^{-(k_{m-1} - k_{m-2})}\\
    & & \cdots (\det Y)^{-k_{1}} + \hbox{terms of lower weight in the $z's$.}
\end{eqnarray*}
where $\kappa = (k_{1}, \dots, k_{m})$.
\end{lem}
\textit{Proof.} If $A$ is a symmetric matrix with latent roots $a_{1}, \dots a_{m}$
then $A^{-1}$ is also a symmetric matrix with latent roots $\alpha_{1}, \dots,
\alpha_{m}$, such that $\alpha_{i} = a_{i}^{-1}$ $i = 1, \dots, m$. Then by
\citet[without proof]{c:66} and \citet[Lemma 2, p.54, with proof]{t:84},
$$
  (\det A)^{n} \frac{C_{\kappa}(A^{-1})}{C_{\kappa}(I_{m})}=\frac{C_{\kappa*}(A)}{C_{\kappa*}(I_{m})}
$$
where $n$ is any integer $\geq k_{1}$ and $\kappa^{*} = (n-k_{m}, \dots n-k_{1})$.
Thus
\begin{eqnarray*}
  C_{\kappa}\left(A^{-1}\right) &=& d_{\kappa}\alpha_{1}^{k_{1}} \cdots \alpha_{m}^{k_{m}}
    + \mbox{ terms of lower weight}\\
    &=& (\det A)^{-n} s_{\kappa,\kappa^{*}}C_{\kappa^{*}}(A) \\
    &=& (\det A)^{-n} s_{\kappa,\kappa^{*}}[d_{\kappa^{*}}a_{1}^{n-k_{m}} \cdots a_{m}^{n-k_{1}}
    + \mbox{ terms of lower weight}]\\
\end{eqnarray*}
Denoting
$$
  s_{\kappa,\kappa^{*}} = \frac{C_{\kappa}(I_{m})}{C_{\kappa^{*}}(I_{m})}
$$
we have
\begin{eqnarray*}
    C_{\kappa}\left(A^{-1}\right) &=& (\det A)^{-n} s_{\kappa,\kappa^{*}}d_{\kappa^{*}}a_{1}^{n-k_{m}-(n -k_{m-1})}
    (a_{1} a_{2})^{n-k_{m-1}-(n-k_{m-2})} \cdots \\
    & & (a_{1}a_{2} \cdots a_{m})^{n-k_{1}} + \cdots \\
    &=& s_{\kappa,\kappa^{*}}d_{\kappa^{*}}a_{1}^{-(k_{m}-k_{m-1})}
    (a_{1} a_{2})^{-(k_{m-1}-k_{m-2})} \cdots (a_{1}a_{2} \cdots a_{m})^{-k_{1}} + \cdots \\
    &=& s_{\kappa,\kappa^{*}}d_{\kappa^{*}}r_{1}^{-(k_{m}-k_{m-1})}
    r_{2}^{-(k_{m-1}-k_{m-2})} \cdots r_{m}^{-k_{1}} + \cdots \\
\end{eqnarray*}
from (39) and (40) ($r_{j} = \tr_{j}(A)$) in \citet[p. 247]{mh:82}
\begin{eqnarray*}
  C_{\kappa}\left(A^{-1}\right) &=& s_{\kappa,\kappa^{*}}d_{\kappa^{*}}\tr_{1}(A)^{-(k_{m}-k_{m-1})}
    \tr_{2}(A)^{-(k_{m-1}-k_{m-2})} \cdots \tr_{m}(A)^{-k_{1}} + \cdots \\
    &=& s_{\kappa,\kappa^{*}}d_{\kappa^{*}} a_{11}^{-(k_{m} - k_{m-1})} \det
  \left[
    \begin{array}{cc}
      a_{11} & a_{12} \\
      a_{21} & a_{22}
    \end{array}
  \right ]^{-(k_{m-1} - k_{m-2})} \cdots (\det A)^{-k_{1}}+ \cdots.
\end{eqnarray*}
Now let $A^{-1} = Y^{-1}Z$, then $A = Z^{-1} Y$ and thus $a_{ij} = z_{i}^{-1}y_{ij}$.
From where
\begin{eqnarray*}
  C_{\kappa}\left(Y^{-1}Z\right) &=& s_{\kappa,\kappa^{*}} d_{\kappa^{*}} (z_{1}^{-1}y_{11})^{-(k_{m} -
  k_{m-1})} \det
  \left[
    \begin{array}{cc}
      z_{1}^{-1}y_{11} & z_{1}^{-1}y_{12} \\
      z_{2}^{-1}y_{21} & z_{2}^{-1}y_{22}
    \end{array}
  \right ]^{-(k_{m-1} - k_{m-2})}\\
    & & \cdots (\det Z^{-1}Y)^{-k_{1}} + \cdots \\
    &=& s_{\kappa,\kappa^{*}}d_{\kappa^{*}}z_{1}^{k_{m}} \cdots z_{m}^{k_{1}} \ y_{11}^{-(k_{m} -
  k_{m-1})} \det
  \left[
    \begin{array}{cc}
      y_{11} & y_{12} \\
      y_{21} & y_{22}
    \end{array}
  \right ]^{-(k_{m-1} - k_{m-2})}\\
    & & \cdots (\det Y)^{-k_{1}} + \cdots,
\end{eqnarray*}
Finally, the result is obtained observe that:
\begin{description}
  \item[i)] In reality, the expression (10) in \citet{c:63}, see also \citet[exprssion (i)(1), p. 228]{mh:82}
     $$
       C_{\kappa}(Y) = d_{\kappa}y_{1}^{k_{1}}...y_{m}^{k_{m}}+\mbox{ terms of lower weight}
     $$
     is
     $$
        C_{\kappa}(Y) = d_{\kappa}(y_{1}^{k_{1}}...y_{m}^{k_{m}}+...+\mbox{
        symmetric terms})+\mbox{ terms of lower weight} \\
     $$
     where ``symmetric terms'' denotes the generic term
     $y_{i_{1}}^{k_{1}}...y_{i_{m}}^{k_{m}}$, with $(i_{1}, \dots, i_{m})$ is a
     permutation of the $m$ integers $1, \dots, m$. Then, alternatively, for a fixed permutation $(i_{1}, \dots,
     i_{m})$,
     $$
       C_{\kappa}(Y) = d_{\kappa}y_{i_{1}}^{k_{1}}...y_{i_{m}}^{k_{m}}+\mbox{ terms of lower weight}
     $$
     or in particular
     $$
       C_{\kappa}(Y) = d_{\kappa}y_{1}^{k_{m}} \cdots, y_{m}^{k_{1}}+\mbox{ terms of lower weight}
     $$

  \item[ii)] Finally, denoting $s_{\kappa,\kappa^{*}}d_{\kappa^{*}}$ by $d_{\kappa}$ and observing that it
     denotes the constant of ``$z_{1}^{k_{1}}...z_{m}^{k_{m}}+...+$ symmetric terms'' of
     $C_{\kappa}(A^{-1})$ $(C_{\kappa}(Y^{-1}Z))$ in terms of the latent roots of $A$ $(YZ^{-1})$. \qed
\end{description}

\begin{rem}
Also, observe that the ``Hint'' in problem 7.5 in \citet{mh:82} is also incorrect.
\end{rem}

An application of Lemma \ref{lemma2}, but in its wrong version, is given by
\citet[Theorem 7.2.13]{mh:82}, surprisingly, the correct result is obtained. The
following results were proposed, without proof, by \citet{c:66} and simultaneously,
with an alternative proof to that given below, by \citet{k:66} and \citet[Lemma 1, p.
53]{t:84}.

\begin{thm} Let $Z$ be a complex symmetric $m \times m$ matrix with $\R(Z) > 0$. Then
\begin{equation}\label{int}
    \int_{X >0} \etr(-XZ) (\det X)^{a -(m+1)/2} C_{\kappa}\left(X^{-1} \right)(dX)
    \hspace{4cm}
\end{equation}
$$\hspace{5cm}
  = \frac{(-1)^{k} \Gamma_{m}[a]}{(-a+(m+1)/2)_{\kappa}} (\det Z)^{-a} C_{\kappa}(Z)
$$
$$\hspace{5.5cm}
  = \frac{ \Gamma_{m}[a]}{(-a+(m+1)/2)_{\kappa}} (\det Z)^{-a} C_{\kappa}(-Z)
$$
for $\R(a) > k_{1} + (m-1)/2$, where $\kappa = (k_{1}, \dots, k_{m})$ and $k = k_{1}
+ \cdots + k_{m}$.
\end{thm}
\textit{Proof.} First suppose that $Z > 0$ is real. Let $f(Z)$ denote the integral on
the left side of (\ref{int}) and make the change of variable $X = Z^{-1/2}YZ^{-1/2}$,
with Jacobian $(dX) = (\det Z)^{-(m+1)/2}(dY)$, to give
\begin{equation}\label{eqint}
    f(Z) = \int_{Y>0} \etr(-Y) (\det Y)^{a-(m+1)/2} C_{\kappa}\left(Y^{-1}Z
    \right)(dY)(\det Z)^{-a}.
\end{equation}
Then, exactly as in the proof of Theorem 7.2.7 in \citet[p. 256-257]{mh:82}
$$
  f(Z) = \frac{f(I_{m})}{C_{\kappa}(I_{m})} C_{\kappa}(Z)(\det Z)^{-a}.
$$
Assuming without loss of generality that $Z = \diag(z_{1},\dots z_{m})$, it then
follows, using (i) from Definition 7.2.1 in \citet{mh:82}, that
\begin{equation}\label{eq1teo}
    f(Z) = \frac{f(I_{m})}{C_{\kappa}(I_{m})}(\det Z)^{-a} d_{\kappa} z_{1}^{k_{1}}\cdots z_{m}^{k_{m}} +
  \mbox{ terms of lower weight}.
\end{equation}
On the other hand, using the result of Lemma \ref{lemma2} in (\ref{eqint}) gives
\begin{eqnarray*}
  f(Z) &=& (\det Z)^{-a} d_{\kappa} z_{1}^{k_{1}}\cdots z_{m}^{k_{m}} \int_{Y>0} \etr(-Y) (\det Y)^{a-(m+1)/2}\\
    & & \times y_{11}^{-(k_{m} -
  k_{m-1})} \det
  \left[
    \begin{array}{cc}
      y_{11} & y_{12} \\
      y_{21} & y_{22}
    \end{array}
  \right ]^{-(k_{m-1} - k_{m-2})} \cdots \det Y^{-k_{1}}(dY) \\
    && + \hbox{ terms of lower weight.}
\end{eqnarray*}
To evaluate this last integral, set $Y = T'T$ where $T$ is upper-triangular with
positive diagonal elements;
$$
  \tr Y = \sum_{i \leq j}^{m}t_{ij}^{2}, \quad y_{11} = t_{ii}^{2} \quad
  \det
  \left[
    \begin{array}{cc}
      y_{11} & y_{12} \\
      y_{21} & y_{22}
    \end{array}
  \right ] = t_{11}^{2}t_{22}^{2},
  \quad \cdots \det Y = \prod_{i=1}^{m}t_{ii}^{2},
$$
and, from Theorem 2.1.9 in \citet[p. 60]{mh:82}
$$
  (dY) = 2^{m}\prod_{i = 1}^{m}t_{ii}^{m+1-i} \bigwedge_{i \leq j}^{m}dt_{ij} = \prod_{i = 1}^{m}\left
  (t_{ii}^{2}\right)^{(m-i)/2} \bigwedge_{i = 1}^{m}dt_{ii}^{2}\bigwedge_{i < j}^{m}dt_{ij}.
$$
Hence
\begin{eqnarray*}
  f(Z) &=& \displaystyle(\det Z)^{-a} d_{\kappa} z_{1}^{k_{1}}\cdots z_{m}^{k_{m}}\prod_{i < j}^{m}
        \left(\int_{-\infty}^{\infty} \exp\left(-t_{ij}^{2}\right)dt_{ij}\right ) \\
    & & \times \prod_{i = 1}^{m} \left (\int_{0}^{\infty} \exp\left(-t_{ii}^{2}\right)
        \left (t_{ii}^{2} \right )^{a-k_{m+1-i}-(i-1)/2- 1} dt_{ii}^{2}
        \right) + \cdots \\
    &=& (\det Z)^{-a} d_{\kappa} z_{1}^{k_{1}}\cdots z_{m}^{k_{m}}
        \pi^{m(m-1)/4}\prod_{i = 1}^{m} \Gamma[a-k_{m+1-i}-(i-1)/2] + \cdots
\end{eqnarray*}
by (\ref{prod3}), we have
$$
   f(Z) = \displaystyle(\det Z)^{-a} d_{\kappa} z_{1}^{k_{1}}\cdots z_{m}^{k_{m}}
            \pi^{m(m-1)/4}\prod_{i = 1}^{m} \Gamma[a-k_{i}-(m-i)/2] + \cdots
$$
which is the result obtained by \citet[eq. (12)]{k:66}, see also \citet[Lemma 1, p.
53]{t:84}. Finally, by (\ref{negpch}) and (\ref{gama2})
\begin{eqnarray}
   f(Z) &=& (\det Z)^{-a} d_{\kappa} z_{1}^{k_{1}}\cdots z_{m}^{k_{m}}
            \pi^{m(m-1)/4}\prod_{i = 1}^{m} \frac{(-1)^{k_{i}}\Gamma[a-(m-i)/2]}{(-a+(m-i)/2+1)_{k_{i}}} + \cdots  \nonumber \\
        &=& (\det Z)^{-a} d_{\kappa} z_{1}^{k_{1}}\cdots z_{m}^{k_{m}}
             \frac{(-1)^{k}\pi^{m(m-1)/4}\prod_{i = 1}^{m}\Gamma[a-(m-i)/2]}{\prod_{i = 1}^{m}(-a+(m-i)/2+1)_{k_{i}}} +
             \cdots \nonumber \\
        &=& (\det Z)^{-a} d_{\kappa} z_{1}^{k_{1}}\cdots z_{m}^{k_{m}}
             \frac{(-1)^{k}\Gamma_{m}[a]}{\prod_{i = 1}^{m}(-a+(m+1)/2-(i-1)/2)_{k_{i}}} + \cdots \nonumber \\
        &=& \label{eq2teo}(\det Z)^{-a} d_{\kappa} z_{1}^{k_{1}}\cdots z_{m}^{k_{m}}
             \frac{(-1)^{k}\Gamma_{m}[a]}{(-a+(m+1)/2)_{\kappa}} + \cdots
\end{eqnarray}
where $k = k_{1}+\cdots + k_{m}$ and $(b)_{\kappa} = \prod_{i=1}^{m}(b - (i -
1)/2)_{k_{i}}$.  Equating coefficients of $z_{1}^{k_{1}}\cdots z_{m}^{k_{m}}$ in
(\ref{eq1teo}) and (\ref{eq2teo}) it follows that
$$
  \frac{f(I_{m})}{C_{\kappa}(I_{m})} =
  \frac{(-1)^{k}\Gamma_{m}[a]}{(-a+(m+1)/2)_{\kappa}}.
$$
Hence we obtain the desired result for real $Z > 0$, and it follows for complex $Z$
with $\R(Z) > 0$ by analytic continuation and recalling that $(-1)^{k}C_{\kappa}(A) =
C_{\kappa}(-A)$. \qed

\begin{rem}
Observe that \citet[penultimate line, p.257]{mh:82} obtains
\begin{eqnarray*}
   f(Z) &=& \displaystyle(\det Z)^{-a} d_{\kappa} z_{1}^{k_{1}}\cdots z_{m}^{k_{m}}
            \pi^{m(m-1)/4}\prod_{i = 1}^{m} \Gamma[a-k_{i}-(i-1)/2] + \cdots\\
        &=& (\det Z)^{-a} d_{\kappa} z_{1}^{k_{1}}\cdots z_{m}^{k_{m}}
            \pi^{m(m-1)/4}\prod_{i = 1}^{m} \frac{(-1)^{k_{i}}\Gamma[a-(i-1)/2]}{(-a+(i-1)/2+1)_{k_{i}}} + \cdots\\
        &=& (\det Z)^{-a} d_{\kappa} z_{1}^{k_{1}}\cdots z_{m}^{k_{m}}
             \frac{(-1)^{k}\pi^{m(m-1)/4}\prod_{i = 1}^{m}\Gamma[a-(m-i)/2]}{\prod_{i = 1}^{m}(-a+(i+1)/2)_{k_{i}}} +
             \cdots\\
        &=& (\det Z)^{-a} d_{\kappa} z_{1}^{k_{1}}\cdots z_{m}^{k_{m}}
             \frac{(-1)^{k}\Gamma_{m}[a]}{\prod_{i = 1}^{m}(-a+(i+1)/2)_{k_{i}}} + \cdots\\
        &\neq& (\det Z)^{-a} d_{\kappa} z_{1}^{k_{1}}\cdots z_{m}^{k_{m}}
             \frac{(-1)^{k}\Gamma_{m}[a]}{(-a+(m+1)/2)_{\kappa}} + \cdots
\end{eqnarray*}
\end{rem}
\begin{cor}Let $V$ be a complex symmetric $m \times m$ matrix with $\R(V) > 0$, and
let $T$ be an arbitrary complex symmetric matrix. Then
\begin{equation}\label{int2}
    \int_{X >0} \etr(-XV) (\det X)^{a -(m+1)/2} C_{\kappa}\left(TX^{-1} \right)(dX)
    \hspace{4cm}
\end{equation}
$$\hspace{5cm}
  = \frac{(-1)^{k} \Gamma_{m}[a]}{(-a+(m+1)/2)_{\kappa}} (\det V)^{-a} C_{\kappa}(VT)
$$
$$\hspace{5.5cm}
  = \frac{ \Gamma_{m}[a]}{(-a+(m+1)/2)_{\kappa}} (\det V)^{-a} C_{\kappa}(-VT)
$$
for $\R(a) > k_{1} + (m-1)/2$, where $\kappa = (k_{1}, \dots, k_{m})$ and $k = k_{1}
+ \cdots + k_{m}$.
\end{cor}
\textit{Proof.} Observe that if $V = I_{m}$ in (\ref{int2}) we obtain (\ref{eqint}).
For the general case substitute $V^{1/2}X V^{1/2}$ for $X$ in (\ref{eqint}) with the
Jacobian of the transformation $|V|^{(m+1)/2}$. \qed

\section*{\large{Conclusions}}
Let us stress that the aim of the present study is not to disparage the importance of
Muirhead's book, but rather to correct the minimal deficiencies we believe to have
identified, and thus help prevent, or minimize, erroneous conclusions being drawn on
the basis of this text, in both current and future work.

\section*{Acknowledgments}
This work was partially supported  by CONACYT-Mexico, research grant 81512 and by
IDI-Spain, grant MTM2005-09209. This paper was written during J. A. D\'{\i}az-
Garc\'{\i}a's stay as a visiting professor at the Department of Statistics and O. R.
of the University of Granada, Spain.

\end{document}